
\input amstex
\documentstyle{amsppt}
\magnification=1200
\overfullrule=0pt
\NoRunningHeads
\document

\centerline{\bf QUANTUM COHOMOLOGY OF A PRODUCT}

\bigskip

\centerline{M. Kontsevich ${}^{1)}$, Yu. Manin ${}^{2)}$}

\smallskip

\centerline{(with Appendix by R. Kaufmann ${}^{2)}$)}

\smallskip

\centerline{1) Math. Dept., University of California at Berkeley, Berkeley CA,
USA}

\centerline{2) Max--Planck--Institut f\"ur Mathematik, Bonn, Germany}

\bigskip

\centerline {\bf 0. Introduction}

\medskip

{\bf 0.1. Quantum cohomology.} Quantum cohomology of a projective algebraic
manifold $V$ is a formal deformation of its cohomology ring. Parameters of this
deformation are coordinates on the space
$H^*(V)$,
and the structure constants are the third derivatives of a formal series
$\Phi^V$
(potential, or free energy) whose coefficients count the number of parametrized
rational curves on $V$ with appropriate incidence conditions (see [KM] for
details).

\smallskip

A natural problem arises how to calculate
$\Phi^{V\times W}$
in terms of
$\Phi^V$
and
$\Phi^W$.
In [KM] it was suggested that this operation corresponds to that of tensor
multiplication of Cohomological Field Theories, or equivalently, algebras over
the moduli operad
$\{ H_* (\overline{M}_{0,n+1})\}$.
The definition of the tensor product depends on a theorem on the structure of
$H_*(\overline{M}_{0,n+1})$
whose proof was only sketched in [KM] (Theorem 7.3). One of the main goals of
this note is to present this proof and related calculations in full detail.
(For another proof, see [G]).

\smallskip

We also discuss the rank one CohFT's and the respective twisting operation. An
interesting geometric example of such a theory is furnished by Weil--Petersson
forms. Potential of this theory is a characteristic function involving
Weil--Petersson volumes calculated in [Z]. We show that a generalization of
WP--forms allows one to construct a canonical coordinate system
on the group of invertible CohFT's.

\smallskip

We start with a brief review of the relevant structures from [KM]. Let $H$ be a
finite-dimensional
$\bold{Z}_2$--graded linear space over a field $K$ of characteristic zero,
endowed with a non--degenerate even symmetric scalar product $g$. The main fact
is the equivalence of two notions:
\roster
\item"i)" A formal solution
$\Phi$
of associativity, or WDVV, equations on
$(H,g)$.
\item"ii)" A structure of CohFT on
$(H,g).$
\endroster

\medskip

{\bf 0.2. {Associativity equations}.} Let
$\{\Delta_a\}$
be a basis of $H$,
$g_{ab} =g(\Delta_a,\Delta_b),\ $
$(g^{ab}) = (g_{ab})^{-1}$.
Denote by
$\gamma =\sum_a x^a \Delta_a$
a generic element of $H$, where
$x^a$
is a formal variable of the same
$\bold{Z}_2$--degree as
$\Delta_a$.
Put
$\partial_a = \partial/\partial x^a$.
A formal series
$\Phi \in K[[x^a]]$
is called a solution of the associativity equations on
$(H,g)$
iff for all
$a,b,c,d$

\TagsOnRight

$$\sum_{ef} \partial_a \partial_b \partial_e \Phi \cdot g^{ef} \partial_f
\partial_c \partial_d \Phi
=(-1)^{a(b+c)} \sum_{ef} \partial_b \partial_c \partial_e \Phi \cdot g^{ef}
\partial_f \partial_a\partial_d \Phi. \tag 0.1$$
Here we use the simplified notation
$(-1)^{a(b+c)}$
for
$(-1)^{\tilde x_a(\tilde x_b +\tilde x_c)}$
where
$\tilde x$
is the
$\bold{Z}_2$--degree of $x$.

\smallskip

We usually assume that
$\Phi$
starts with terms of degree
$\ge 3$,
or identify
$\Phi$
and
$\Phi^{\prime}$
differing by a polynomial of degree
$\le 2$.
An extensive geometric treatment of the associativity equations is given in
[D]. For the next definition, remind that
$\overline M_{0n}$
denotes the moduli space of stable curves of genus zero with $n$ labelled
pairwise distinct points: see [Ke].

\medskip

{\bf 0.3. Cohomological field theories.} A structure of the tree level
Cohomological Field Theory on
$(H,g)$
is given by a sequence of
$S_n$--covariant $K$--linear maps
$$I_n : H^{\otimes n} \rightarrow H^* (\overline{M}_{0n},K),\,\  n\ge 3 \tag
0.2$$
satisfying the following set of identities (0.3). The values of
$I_n$
generally are not homogeneous.

\smallskip

Consider an unordered partition
$\sigma :\ \{1,\ldots,n\} =S_1 \coprod S_2,\, |S_i|\ge 2$.
It defines an embedding of the boundary divisor
$\varphi_{\sigma} : \overline{M}_{0,n_1+1} \times \overline{M}_{0,n_2+1}
\rightarrow \overline{M}_{0,n}$.
Over the generic point of this divisor the universal curve consists of two
components, and the labelled points are distributed between them according to
$\sigma$.
The maps (0.2) must satisfy for all
$n\ge 3$
the relations:
$$\varphi^*_{\sigma} (I_n (\gamma_1 \otimes \ldots \otimes \gamma_n)) =\epsilon
(\sigma) (I_{n_1+1} \otimes I_{n_2+1}) (\otimes_{j\in S_1} \gamma_j \otimes
\Delta \otimes (\otimes_{k \in S_2} \gamma_k)) \tag 0.3$$
where
$\Delta =\Sigma \Delta_a \otimes \Delta_b g^{ab}$
is the Casimir element, and
$\epsilon (\sigma)$
is the sign of the permutation induced on the odd arguments
$\gamma_1, \ldots, \gamma_n$.
There are two other useful reformulations of CohFT. First, dualizing (0.2) we
get a series of maps

$$I^t_{n+1} : H_* (\overline{M}_{0,n+1},K) \rightarrow Hom(H^{\otimes n}, H),\,
\ n\ge 2. \tag 0.4$$
Thus any homology class in
$\overline{M}_{0,n+1}$
is interpreted as an $n$-ary opertaion on $H$. The relations (0.3) become
identities between these operations whose totality means that $H$ is given a
structure of an algebra over the cyclic operad
$\{H_*(\overline{M}_{0,n+1},K)\}$
(see [GK]).

\smallskip

Second, we can iterate the maps
$\varphi_{\sigma}$
in order to study the restrictions of the classes
$I_n(\gamma_1 \otimes \ldots \otimes \gamma_n)$
to all boundary strata. These strata are naturally indexed by the (dual) trees
of stable curves, which form a category. Both sides of (0.2) extend to the
functors on this category, and (0.3) says that
$\{I_n\}$
becomes a functor morphism. Below we will extensively use the combinatorial
side of this picture explained in [KM].

\medskip

{\bf 0.4. From $\{I_n\}$ to $\Phi$.} Every CohFT
$\{H,g,I_n\}$
defines a sequence of symmetric polynomials
$Y_n : M^{\otimes n} \rightarrow K$:

$$Y_n (\gamma_1 \otimes \ldots \otimes \gamma_n) := \int_{\overline{M}_{0n}}
I_n (\gamma_1 \otimes \ldots \otimes \gamma_n) \tag 0.5$$
Put $x=(x^a)$ and

$$\Phi (x) =\sum_{n\ge 3} \frac 1 {n!} Y_n ((\sum_a x^a \Delta_a)^{\otimes n}).
\tag 0.6$$
 From (0.3) and Keel's linear relations between boundary divisors in
$\overline{M}_{0n}$
one can formally deduce (0.1).

\medskip

{\bf 0.5. From $\Phi$ to $\{I_n\}$.}
This transition is markedly more difficult. If
$I_n$
with the properties (0.5), (0.6) exist at all, they are defined uniquely,
because iterating (0.3) one can calculate integrals of
$I_n$
over all boundary strata of
$\overline{M}_{0n}$.
Namely, let
$\tau$
be the dual tree of a stable curve $C$ of genus zero with
$n$ marked points. Let us remind that the set of vertices
$V_{\tau}$
consists of irreducible components of $C$, edges
$E_{\tau}$
are (in a bijection with) double points of $C$, tails
$T_{\tau}$
(one vertex edges) are marked points. Incidence relations between
$V_{\tau}, E_{\tau},T_{\tau}$
reflect those in $C$. A flag of
$\tau$
is a pair (vertex, incident edge or tail); the set of flags is denoted
$F_{\tau}$.
Let
$\overline{M}_{\tau} \subset \overline{M}_{0n}$
be the submanifold parametrizing curves of type
$\tau$
and their specialisations. We have
$\overline{M}_{\tau} \cong \prod_{v \in V_{\tau}} \overline{M}_{|v|}$.
Here $|v|=|F_{\tau}(v)|,\, F_{\tau}(v)$ is the set of flags
incident to $v.$
The homology classes of all
$\overline{M}_{\tau}$
generate
$H_*(\overline{M}_{0n},K)$,
and we have from (0.3):

$$\int_{\overline{M}_{\tau}} I_n (\gamma_1 \otimes \ldots \otimes \gamma_n) =
(\bigotimes_{v \in V_{\tau}} Y_{|v|}) (\gamma_1 \otimes  \ldots \gamma_n
\otimes \Delta ^{\otimes E_{\tau}}) \in K .\tag 0.7$$
We use here the formalism of tensor products indexed by arbitrary finite sets
and interpret the argument of the r.h.s. (0.7) as an element of
$H^{\otimes F_{\tau}}$
and
$\bigotimes_{v\in V_{\tau}} Y_{|v|}$
as a function on
$H^{\otimes F_{\tau}}$.
(0.3) is a particular case of (0.7) for an one-edge tree.

\smallskip

In order to use (0.7) for construction of
$\{I_n\}$
it remains to check that the r.h.s. of (0.7) satisfies all the linear relations
between the classes of
$\overline{M}_{\tau}$
in
$H_*(\overline{M}_{0n})$.
This was made in sec. 8 of [KM] modulo the theorem, describing these relations
and proved in \S \, 2 of the present paper.

\medskip

{\bf 0.6. The tensor product.}
Let
$\{H^{\prime}, g^{\prime}, I^{\prime}_n\}$
and
$\{H^{\prime\prime},g^{\prime\prime},I^{\prime\prime}_n\}$
be two CohFT's. Put
$H = H^{\prime} \otimes H^{\prime\prime}$
and
$g =g^{\prime} \otimes g^{\prime\prime}$.
We can define a CohFT on
$(H,g)$
by
$$I_n(\gamma^{\prime}_1 \otimes \gamma^{\prime\prime}_1 \otimes \ldots \otimes
\gamma^{\prime}_n \otimes \gamma^{\prime\prime}_n) := \epsilon
(\gamma^{\prime}, \gamma^{\prime\prime})I^{\prime}_n (\gamma^{\prime}_1 \otimes
\ldots \otimes \gamma^{\prime}_n) \wedge I^{\prime\prime}_n
(\gamma^{\prime\prime}_1 \otimes \ldots \otimes \gamma^{\prime\prime}_n) \tag
0.8$$
where
$\epsilon (\gamma^{\prime},\gamma^{\prime\prime})$
is the standart sign in superalgebra, and
$\wedge$
is the cup product in
$H^*(\overline{M}_{0n},K)$.
One can easily check (0.3). Although (0.8) looks very simple on the level of
full CohFT's, it cannot be trivially restricted to calculate the potential. In
fact, for potential we must know only the higher dimensional component of
$I_n$
(see (0.5)), but it involves components of all degrees of
$I^{\prime}_n,I^{\prime\prime}_n$
(see (0.8)), which are given by (0.7) as functionals on the boundary homology
classes.

\smallskip

The only remaining obstruction to calculating
$\Phi$
of
$\{I_n\}$
is thus our incomplete understanding of the cup product in terms of dual
boundary classes. Ralph Kaufmann obtained a fairly simple formula for the
intersection indices of strata of complementary dimensions. With his
permission, we reproduce it in \S \, 2. But at the moment we are unable to
invert this Gram matrix (of redundant size).

\medskip

{\bf 0.7. Application to quantum cohomology.}
In the application of this formalism to the quantum cohomology of $V$ we have:
$H = H^*(V,K),\, g=\ $
Poincar\'e pairing, and
$\int_{\overline{M}_{0n}} I_n (\gamma_1 \otimes \ldots \otimes \gamma_n)$
is the (appropriately defined) number of ``stable maps''
$(C,x_1,\ldots,x_n;\varphi)$
where $C$ is a curve of genus 0 with $n$ marked points
$x_1,\ldots,x_n,\, \varphi : C \rightarrow V$
is a morphism, such that any connected component contracted by
$\varphi$
(together with its special points) is stable. These Gromov-Witten invariants
are of primary interest in enumerative geometry, and we expect that the tensor
product of CohFT's described above furnishes an algorithm for calculating them
on
$V \times W$
from those on $V$ and $W$. Strictly speaking, this must be proved starting with
a geometric construction of GW--invariants. In concrete examples it suffices to
simply check the coincidence of a finite set of coefficients of the two
potentials using The First Reconstruction Theorem 3.1 of [KM]. For example,
$\Phi^{\bold{P}^1 \times \bold{P}^1}$
actually coincides with the tensor square of the simple potential
$\Phi^{\bold{P}^1}$.
Namely, if
$\Delta_0 \in H^0(\bold{P}^1)$
and
$\Delta_1 \in H^1(\bold{P}^1)$
are respectively the fundamental class and the dual class of a point, we have

$$\alignat 1 &\Phi^{\bold{P}^1} (x\Delta_0 +z\Delta_1) =\frac 1 2 x^2z
+e^z-1-z-\frac {z^2}2 , \tag 0.8\\
\text{  }\\
&\Phi^{\bold{P}^1\times \bold{P}^1} (x\Delta^{\otimes 2}_0 +y^1\Delta_1 \otimes
\Delta_0 +y^2 \Delta_0 \otimes \Delta_1 + z\Delta_1^{\otimes 2}) = \\
&= \frac 1 2 x^2z +xy^1y^2 + \sum_{a+b \ge 1} N(a,b) \frac {z^{2a+2b-1}}
{(2a+2b-1)!} \, e^{ay^1+by^2} \tag 0.9
\endalignat$$
where
$N(a,b)$
is the number of rational curves of bidegree
$(a,b)$
on
$\bold{P}^1 \times \bold{P}^1$
passing through
$2a+2b-1$
points in general position. As is remarked in [DFI], the structure of
$N(a,b)$
which can be derived from the associativity equations looks simpler than that
of GW--numbers for
$\bold{P}^2$,
probably because of this tensor product property.

\smallskip

In the remaining part of this Introduction, we summarize our notation and
conventions about the combinatorics of trees and the (co)homology of
$\overline{M}_{0n}$.
The next section is devoted to the multiplicative properties of strata classes.
In \S \, 2 we prove the completeness of the standard linear relations between
them. The last \S \, 3 discusses rank one CohFT's.

\medskip

{\bf 0.8. Partitions and trees.}
As in [KM], a tree
$\tau$
for us is a system of finite sets
$(V_{\tau},E_{\tau},T_{\tau}),\, V_{\tau} \ne \emptyset$,
with appropriate incidence relations defining
$F_{\tau}$
(see 0.5 above for notation). A structure of $S$-tree on
$\tau$
(where $S$ is a finite set) is given by a bijection
$T_{\tau} \rightarrow S$.
Sometimes we identify
$T_{\tau}$
with $S$ using this bijection. A tree
$\tau$
is {\it stable} if
$|v| \ge 3$
for all
$v \in V_{\tau}$.
Most our trees are stable.

\smallskip

A stable $S$-tree
$\tau$
corresponds to a (family of) stable curve(s) of genus 0 with points labelled by
elements of $S$. One-edge $S$-trees are in a bijection with unordered
partitions of $S$ into two subsets
$\sigma :\ S =S_1 \coprod S_2$;
stability means that
$|S_i| \ge 2$;
the tails marked by
$S_i$
belong to the vertex
$v_i,\, i=1,2$.
We will systematically identify such partitions with the corresponding trees.
For such a
$\sigma$
and, say, four elements
$i,j,k,l \in S$
we use a notation like
$ij\sigma kl$
to imply that
$\{i,j\}$
and
$\{k,l\}$
belong to different parts of
$\sigma$.

\smallskip

For two unordered stable partitions
$\sigma =\{S_1,S_2\}$
and
$\tau =\{T_1,T_2\}$
of $S$ put
$$\align a(\sigma,\tau):= \quad &\text{the number of non-empty pairwise} \\
&\text{distinct sets among} \quad S_i \cap T_j,\, i,j =1,2. \endalign$$
Clearly,
$a(\sigma,\tau)=2,3$,
or $4$. Moreover,
$a(\sigma,\tau)=2$
iff
$\sigma =\tau$,
and
$a(\sigma,\tau) =4$
iff there exist pairwise distinct
$i,j,k,l \in S$
such that simultaneously
$ij \sigma kl$
and
$ik \tau jl$.
If
$a(\sigma,\tau) =3$,
we sometimes call
$\sigma$
and
$\tau$
compatible. A family of 2-partitions
$\{\sigma_1,\ldots,\sigma_m\}$
is called {\it good}, if for all
$i \ne j,\, \sigma_i$
and
$\sigma_j$
are compatible. $S$-trees form objects of several catgeories differing by the
size of their morphism sets. The most useful morphisms
$f : \tau \rightarrow \sigma$
contract several edges and tails of
$\tau$:
$f$ induces a surjection
$V_{\tau} \rightarrow V_{\sigma}$
and injections
$E_{\sigma} \rightarrow E_{\tau},\, T_{\sigma} \rightarrow T_{\tau}$;
labelling sets for
$\tau$
and
$\sigma$
may differ. We will mostly consider morphisms of $S$-trees identical on $S$
(pure contractions of edges, or $S$-moprhisms). The one-vertex tree is a final
object in the category of $S$-trees and $S$-morphisms. If a direct product
$\sigma \times \tau$
of two $S$-trees in this category exists, it comes equipped with two
contractions
$\sigma \times \tau \rightarrow \sigma$
and
$\sigma \times \tau \rightarrow \tau$.
A geometrically nice case is when
$|E_{\sigma \times \tau}| =|E_{\sigma}| + |E_{\tau}|$.
E.g. for one-edge trees this is the case when
$a(\sigma,\tau) =3$.
For
$a(\sigma,\tau) =4,\, \sigma \times \tau$
does not exist, and for
$\sigma =\tau$
we have
$\sigma \times \sigma =\sigma$.

\smallskip

A few more words about the geometry of an individual tree
$\tau$. Any flag
$f =(v,e),\, v \in V_{\tau}, \, e \in E_{\tau}$
or $e \in T_{\tau}$
defines a complete subgraph
$\beta (f)$
of
$\tau$
which we will call {\it the branch} of $f$. If
$e \in T_{\tau},\, \beta (f)$
consist of the vertex $v$ and tail $e$. Generally,
$\beta (f)$
includes
$v,e$,
and all edges, vertices and tails that can be reached from $v$ by a no-return
path starting with $f$. We denote by
$T_{\tau} (f)$
the tails belonging to
$\beta (f)$,
and by
$S(f)$
their labels (if
$\tau$
is a labelled tree).

\medskip

{\bf 0.9. Moduli spaces.}
For a finite set
$S,\, |S| \ge 3,\, \overline{M}_{0S}$
parametrizes stable curves of genus zero with a family of pairwise distinct
points labelled by $S$. More generally, for a stable $S$-tree
$\tau,\, \overline{M}_{\tau}$
parametrizes such curves with dual graph (isomorphic to)
$\tau;\,  \overline{M}_{0S}$
corresponds to the one-vertex $S$-tree. Any pure contraction
$\tau \rightarrow \sigma$
bijective (but not necessarily identical) on $S$ induces a morphism
$\overline{M}_{\tau} \rightarrow \overline{M}_{\sigma}$.
In this way
$\{\overline{M}_{\sigma}\}$
form a topological cyclic operad (see [GK]), and
$\{H_*( \overline{M}_{0S})\}$
form a linear cyclic operad.

\smallskip

If an $S$-morphism of $S$-trees
$\tau \rightarrow \sigma$
exists, it is unique, and
$\overline{M}_{\tau}\rightarrow \overline{M}_{\sigma}$
is a closed embedding whose image is called a (closed) stratum of
$\overline{M}_{\sigma}$.
In particular, all
$\overline{M}_{\tau}$
``are'' closed strata in
$\overline{M}_{0S}$.
We have
$\overline{M}_{\tau} \cong \prod_{v\in V_{\tau}} \overline{M}_{0F_{\tau} (v)}$
(canonically),
$\overline{M}_{0,F_{\tau}(v)} \cong \overline{M}_{0|v|}$
(non-canonically). The codimension of the stratum
$\overline{M}_{\tau}$
in
$\overline{M}_{\tau}$
is
$|E_{\tau}|$.
In particular, stable one-edge $S$-trees (and stable 2-partitions)
$\sigma$
bijectively correspond to the boundary divisors.

\medskip

{\bf 0.10. Keel's presentation.}
Fixing
$S,\, |S| \ge 3$,
we denote by
$\{D_{\sigma}|\sigma$
stable 2-partitions of
$S\}$
a family of commuting independent variables. Put
$F_S = K[D_{\sigma}]\ (F_S=K$
for
$|S|=3)$.
We consider
$F_S$
as a graded polynomial ring,
$\roman{deg} \, D_{\sigma} =1$.
Define the ideal
$I_S \subset F_S$
by means of the following generators:
\roster
\item"a)" For each pairwise distinct foursome
$i,j,k,l \in S$:

$$R_{ijkl} := \sum_{ij\sigma kl} D_{\sigma} - \sum_{kj\tau il} D_{\tau} \in I_S
. \tag 0.10$$
\item"b)" For each pair
$\sigma,\tau$
with
$a(\sigma,\tau) =4$:

$$D_{\sigma} D_{\tau} \in I_S .\tag 0.11$$
\endroster
Finally, put
$H^*_S =K[D_{\sigma}]/I_S$.

\medskip

\proclaim {\quad 0.10.1. Theorem (Keel [Ke])}
The map
$$\align D_{\sigma} \longmapsto \quad &\text{dual cohomology class of the
boundary divisor} \\
&\text{in $\overline{M}_{0S}$ corresponding to the partition $\sigma$}
\endalign$$
induces the isomorphism of rings (doubling the degrees)

$$H^*_S \overset \sim \to \longrightarrow H^* (\overline{M}_{0S},K). \tag
0.12$$
\endproclaim

\smallskip

Since
$\overline{M}_{0S}$
is a smooth manifold whose homology and cohomology is generated by algebraic
classes, on which homological and rational equivalences coincide, (0.12)
describes the homology and the Chow ring as well. In addition, Keel's
presentation is very convenient for describing the operadic structure maps.
E.g. bijections
$S^{\prime} \rightarrow S^{\prime\prime}$
(relabelling of points) translate simply by the respective relabelling of
$D_{\sigma}$'s.
For the remaining morphisms, see the next section.

\newpage

\centerline{\bf \S \, 1. Boundary strata and the multiplicative structure of
$H^*(\overline{M}_{0S})$}

\bigskip

{\bf 1.1. Good monomials.}
The monomial
$D_{\sigma_1} \ldots D_{\sigma_a} \in F_S$
is called good, if the family of 2-partitions
$\{\sigma_1,\ldots,\sigma_a\}$
is good, i.e.
$a(\sigma_i,\sigma_j)=3$
for
$i \ne j$.
In particular,
$D_{\sigma}$
and 1 are good.

\medskip

\proclaim{\quad 1.2. Lemma}
Let
$\tau$
be a stable $S$-tree with
$|E_{\tau}| \ge 1$.
For each
$e \in E_{\tau}$,
denote by
$\sigma (e)$
the 2-partition of $S$ corresponding to the one edge $S$-tree obtained by
contracting all edges except for $e$. Then

$$m(\tau) := \prod_{e\in E_{\tau}} D_{\sigma (e)}$$
is a good monomial.
\endproclaim

\smallskip

{\bf Proof.}
Let
$e \ne e^{\prime} \in E_{\tau}$.
There exists a sequence of pairwise distinct edges
$e=e^{\prime}_0,e^{\prime}_1,\ldots, e^{\prime}_r,\, e^{\prime}_{r+1}
=e^{\prime},\, r\ge 0$,
such that
$e^{\prime}_j$
and
$e^{\prime}_{j+1}$
have a common vertex
$v_j$.

\midspace{4cm}\caption{Figure 1. Arrows symbolize branches}

Let $u$ be the remaining vertex of $e,\, w$
that of
$e^{\prime}$.
Let
$S^{\prime}$
be the set of all tails of
$\tau$
belonging to the branches starting at $u$ but not with a flag belonging to
$e_0$;
similarly, let
$S^{\prime\prime}$
be the set of all tails of
$\tau$
belonging to the branches that start at $w$ but not with a flag belonging to
$e^{\prime}$.
Finally, let $T$ be the set of all tails on the branches at
$v_0,\ldots, v_r$
not starting with the flags in
$e^{\prime}_0,\ldots, e^{\prime}_{r+1}$
(we identify tails with their labels). Since
$\tau$
is stable, all three sets
$S^{\prime}, S^{\prime\prime}$
and $T$ are non-empty. Finally

$$\sigma (e) =\{ S^{\prime}, S^{\prime\prime} \coprod T\},\, \sigma
(e^{\prime}) =\{S^{\prime} \coprod T,\, S^{\prime\prime}\}.$$
It follows that
$a(\sigma (e),\, \sigma (e^{\prime})) =3$
so that
$m(\tau)$
is a good monomial.$\hfill \blacksquare$

\smallskip

We put
$m(\tau) =1$,
if
$|E_{\tau}| =0$.

\medskip

\proclaim{\quad 1.3. Proposition}
For any
$1 \le r \le |S|-3$,
the map
$\tau \longmapsto m(\tau)$
establishes a bijection between the set of good monomials of degree $r$ in
$F_S$
and stable $S$-trees
$\tau$
with
$|E_{\tau}| =r$
modulo $S$-isomorphism. There are no good monomials of degree
$> |S| -3$.
\endproclaim

\smallskip

{\bf Proof.}
For
$r=0,1$
the assertion is clear. Assume that for some
$r \ge 1$
the map
$\tau \longmapsto m(\tau)$
is surjective on good monomials of degree $r$. We will prove then that it is
surjective in the degree
$r+1$.

\smallskip

Let
$\roman{deg} \, m^{\prime} =r+1$.
Choose a divisor
$D_{\sigma}$
of
$m^{\prime}$
which is {\it extremal} in the following sense: one element, say
$S_1$,
of the partition
$\sigma =\{S_1,S_2\}$
is minimal in the set of all elements of all 2-partitions
$\sigma^{\prime}$
such that
$D_{\sigma^{\prime}}$
divides
$m^{\prime}$.
Put
$m^{\prime} =D_{\sigma} m$.
Since $m$ is good of degree $r$, we have
$m =m(\tau)$
for some stable $S$-tree
$\tau$.
We will show that
$m^{\prime} =m(\tau^{\prime})$
where
$\tau^{\prime}$
is obtained from
$\tau$
by inserting a new edge with tails marked by
$S_1$
at an appropriate vertex
$v \in V_{\tau}$.
(This means that vice versa, there exists a contraction
$\tau^{\prime} \rightarrow \tau$
of one edge to the vertex having incident flags
$S_1$
and (half of) this edge).

\midspace{4cm} \caption{Figure 2. Inserting edge at a vertex}

First we must find $v$ in
$\tau$.
To this end, consider any edge
$e \in E_{\tau}$
and the respective partition
$\sigma (e) :\{ S^{\prime}_e, S^{\prime\prime}_e\}$
(obtained by contracting all edges except for $e$). Since
$a(\{S_1,S_2\},\, \{S^{\prime}_e,S^{\prime\prime}_e\}) =3$
and
$S_1$
is minimal, one sees that exactly one of the sets
$\{S^{\prime}_e,S^{\prime\prime}_e\}$
strictly contains
$S_1$.
Let it be
$S^{\prime\prime}_e$.
Orient $e$ by declaring that the direction from the vertex (corresponding to)
$S^{\prime}_e$
to
$S^{\prime\prime}_e$
is positive. We claim that with this orientation, for any
$w \in V_{\tau}$
there can be at most one edge outgoing from $w$. In fact, if
$\tau$
contains a vertex $w$ with two positively oriented flags
$f_1$and $f_2$, then $S_1$ must be contained in the two subsets
of $S,$ $S(f_1)$ and $S(f_2)$. But their intersection is empty.

\midspace{4cm} \caption{Figure 3.}

It follows that there exists exactly one vertex $v\in V_{\tau}$ having no
outgoing edges. Moreover, $S_1$ is contained in the set of labels
of the tails at $v$ by construction. If we now define
$\tau^{\prime}$ by inserting a new edge $e^{\prime}$ at $v$ so
that $\sigma (e^{\prime})=\sigma$, we will clearly have
$m^{\prime}=m(\tau^{\prime}).$ If $r\le |S|-4,$ the tree
$\tau^{\prime}$ cannot be unstable because, first, $|S_1|\ge 2,$
and second, at least two more flags converge at
$|v|$: otherwise the unique incoming edge would produce the partition
$\{ S_1,S_2\}=\sigma$ which would mean that $D_{\sigma}$
divides already $m(\tau ).$

\smallskip

For $r=|S|-3,$ this argument shows that $m^{\prime}$ cannot exist
because all the vertices of $\tau$ have valency three.

\smallskip

It remains to check that if $m(\tau_1)=m(\tau_2),$ then $\tau_1$
and $\tau_2$ are isomorphic.

\smallskip

Assume that this has been checked in degree $\le r$ and that
$\roman{deg}\ \tau_1=\roman{deg}\ \tau_2=r+1.$ Choose an extremal divisor
$D_{\sigma}$ of $m(\tau_1)=m(\tau_2)$ as above and contract
the respective edges of $\tau_1,\tau_2$ getting the trees
 $\tau_1^{\prime},\tau_2^{\prime}.$ Since
$m(\tau_1^{\prime})=m(\tau_2^{\prime})
=m(\tau_i)/D_{\sigma}$, $\tau_1^{\prime}$ and $\tau_2^{\prime}$
are isomorphic by the inductive assumption. This isomorphism respects
the marked vertices $v_1^{\prime},v_2^{\prime}$ corresponding to the
contracted edges because as we have seen they are uniquely defined.
Hence it extends to an $S$-isomorphism
$\tau_1 \rightarrow \tau_2$. $\hfill \blacksquare$

\medskip

{\bf 1.3.1. Remark.}
Proposition 1.3 and Keel's theorem 0.10.1, together with the fact that the
boundary divisors have transversal intersections, show that the image of a good
monomial
$m(\tau)$
in
$H^*(\overline{M}_{0S})$
is the dual class of the stratum
$\overline{M}_{\tau}$.

\medskip

{\bf 1.4. Multiplication formulas I.}
Let now
$\sigma,\tau$
be two stable $S$-trees,
$|E_{\sigma}| =1$.
We have the following three possibilities

\smallskip

a). $D_{\sigma} m(\tau)$
{\it is a good monomial}. Then

$$D_{\sigma} m(\tau) =m(\tau^{\prime}) \tag 1.1$$
where
$\tau^{\prime} \rightarrow \tau$
is the unique $S$-morphism contracting the edge in
$E_{\tau^{\prime}}$,
whose 2-partition coincides with that of
$\sigma$.

\smallskip

More generally, if
$m(\sigma) m(\tau)$
is a good monomial, then
$$m(\sigma) m(\tau) =m(\sigma \times \tau) \tag 1.2$$
where the direct product is the categorical one in the category of $S$-trees
and $S$-morphisms. We can identify
$E_{\sigma \times \tau}$
with
$E_{\sigma} \coprod E_{\tau}$,
and
$p_1:\ \sigma \times \tau \rightarrow \sigma$
(resp.
$p_2:\ \rho \times \tau \rightarrow \tau$)
contracts edges of the second factor (resp. of the first one).

\smallskip

b). {\it There exists a divisor} $D_{\sigma^{\prime}}$ {\it of}
$m(\tau), \, |E_{\sigma^{\prime}}| =1$, {\it such that}
$a(\sigma,\sigma^{\prime}) =4$.
Then

$$D_{\sigma} m(\tau) \equiv 0 \ \roman{mod} \ I_S, \tag 1.3$$
where
$I_S \subset F_S$
is the ideal of Keel's relations.

\smallskip

c). $D_{\sigma}$ {\it divides} $m(\tau)$.
Then let
$e \in E_{\tau}$
be the edge corresponding  to
$\sigma;\, v_1,\, v_2$
its vertices,
$(v_i,e)$
the corresponding flags.

\smallskip

We will write several different expressions for
$D_{\sigma} m(\tau) \ \roman{mod} \ I_S$,
corresponding to various possible choices of unordered pairs of distinct flags
$\{ \bar \iota, \bar j\} \subset F_{\tau} (v_1) \setminus \{(v_1,e)\},\, \{\bar
k, \bar l\} \subset F_{\tau} (v_2) \setminus \{(v_2,e)\}$.
For each choice, put
$$\align &T_1 =F_{\tau}(v_1) \setminus \{\bar \iota, \bar j, (v_1,e)\}, \\
&T_2 = F_{\tau} (v_2) \setminus \{ \bar k, \bar l, (v_2,e)\}. \endalign $$
Notice that because of stability the set of such choices is non-empty.

\smallskip

\proclaim{\quad 1.4.1. Proposition}
For every such choice we have
$$D_{\sigma} m(\tau) \equiv - \sum \Sb T \subset T_1 \\ |T| \ge 1 \endSb
m(tr_{T,e} (\tau)) - \sum \Sb T \subset T_2 \\ |T| \ge 1 \endSb m(tr_{T,e}
(\tau ))\ \roman{mod} \ I_S \tag 1.4$$
where
$tr_{T,e}(\tau)$
is the tree obtained from
$\tau$
by ``transplanting all branches starting in $T$ to the middle point of the edge
$e$.'' (An empty sum is zero).
\endproclaim

\midspace{6cm} \caption{Figure 4. Transplants: arrows symbolize branches}

{\bf Remark.} We can also describe
$tr_{T,e} (\tau)$
as a result of inserting an extra edge instead of the vertex
$v_1$
(resp.
$v_2$) and putting the branches $T$ to the common vertex of the new edge and
$e$, similarly to what we have done in the proof of Proposition 1.3. There
exists a unique $S$-morphism
$tr_{T,e} (\tau) \rightarrow \tau$
contracting  one edge.

\smallskip

{\bf Proof.} We choose pairwise distinct labels on the chosen branches
$i \in S(\bar \iota),\, j \in S(\bar j),\, k \in S(\bar k),\, l \in S(\bar
\ell)$
and then calculate the element (see (0.10))

$$R_{ijkl} \cdot m(\tau) =\left( \sum_{ij\rho kl} D_{\rho} - \sum_{kj\rho il}
D_{\rho} \right)m(\tau) \equiv 0 \ \roman{mod} \ I_S .\tag 1.5$$
Since
$ij\sigma kl,$
for all terms
$D_{\rho}$
of the second sum in (1.5) we have
$a(\sigma,\rho) =4$
so that
$D_{\rho}m(\tau ) \in I_S$.
Among the terms of the first sum, there is one
$D_{\sigma}$.
If
$ij\rho kl$
and
$\rho \ne \sigma$,
then
$D_{\rho}$
cannot divide
$m(\tau).$
Otherwise
$\rho$
would correspond to an edge
$e^{\prime} \ne e$,
but the 2-partition of such an edge cannot break
$\{i,j,k,l\}$
into
$\{i,j\}$
and
$\{k,l\}$
as a glance to a picture of
$\tau$
shows. It follows that
$D_{\rho} m(\tau) = m(\rho \times \tau)$
as in (1.2). The projection
$\rho \times \tau \rightarrow \tau$
contracts the extra edge onto a vertex that can be only one of the ends of $e$,
otherwise, as above, the condition
$ij\rho kl$
cannot hold. It should be clear by now that
$\rho \times \tau$
must be one of the trees
$tr_{T,e} (\tau)$,
and that each tree of this kind can be uniquely represented as
$\rho \times \tau$
for some
$\rho$
with
$ij\rho kl$.
But from (1.5) it follows that

$$D_{\sigma} m(\tau) \equiv -\sum \Sb ij\rho kl \\ \rho \ne \sigma \endSb
D_{\rho} m(\tau) \ \roman{mod} \ I_S $$
which is (1.4). $\hfill \blacksquare$

\medskip

\proclaim{\quad 1.5. Corollary}
Classes of good monomials linearly generate
$F_S/I_S =H^*_S$.

\endproclaim

This follows from (1.1), (1.3) and (1.4) by induction on the degree.

\smallskip

{\bf 1.5.1. Remark.} Formulas (1.1)--(1.4) (and (1.7) below) can be rewritten
as expressing operadic morphisms
$\varphi^*_{\sigma} : H^*(\overline{M}_{0S}) \rightarrow
H^*(\overline{M}_{\sigma}) = H^*(\overline{M}_{0,S_1\coprod \{\cdot \}} )
\otimes H^* (\overline{M}_{0,S_2 \coprod \{\cdot\}})$,
$\sigma =\{S_1,S_2\}$
in terms of classes of boundary divisors. E.g. (1.4) means that

$$\varphi^*_{\sigma} (D_{\sigma}) \equiv -\sum_{ij\tau \{\cdot\}} D_{\tau}
\otimes1 -\sum_{\{\cdot \} \rho \{kl\}} 1 \otimes D_{\rho}\ ,\tag 1.6$$
where $\tau$ (resp. $\rho$) runs over stable 2-partitions of $S_1\coprod \{.\}$
(resp. $S_1\coprod \{.\}$).

\medskip

{\bf 1.6. Multiplication formulas II.} It may be more convenient to have
formulas independent of arbitrary choice of $\{\overline{i}, \overline{j},
\overline{k}, \overline{l}\} .$ One way to achieve this is to average
(1.4) over all possible choices. We will illustrate this procedure by
calculating $D_{\sigma}^2$ and some of the ``tautological classes.''

\smallskip

Let
$\sigma$
be a stable 2-partition
$\{T_1,T_2\}$
of $S$.

\proclaim{\quad 1.6.1. Proposition}
We have

$$\align D^2_{\sigma} = &-\sum \Sb T\subset T_1 \\ 1 \le |T|\le |T_1|- 2 \endSb
D_{\sigma} D_{\{T_1\setminus T,T_2 \coprod T\}} \frac {|T_1\setminus
T|(|T_1\setminus T|-1)} {|T_1|(|T_1|- 1)} \\
&-\sum \Sb T\subset T_2 \\ 1 \le |T|\le |T_2|- 2 \endSb D_{\sigma}
D_{\{T_1\coprod T,T_2 \setminus T\}} \frac {|T_2\setminus T|(|T_2\setminus
T|-1)} {|T_2|(|T_2|- 1)} \ \roman{mod} \ I_S \tag 1.7\endalign$$
\endproclaim

{\bf Proof.} We first write
$\pmatrix |T_1| \\2 \endpmatrix$
$\pmatrix |T_2| \\ 2 \endpmatrix$
identities (1.4) for all possible choices of
$i,j \in T_1,\, k,l \in T_2$,
then sum them up and change the summation order by first choosing subsets
$T \subset T_1$
or
$T_2$,
and then
$i,j$
or
$k,l$
in the complement. $\hfill \blacksquare$

\medskip

{\bf 1.6.2. Tautological classes.} These classes
$\tau_d^{(i)} \in H^*(\overline{M}_{0n})$
are defined as
$c_1(T^*_{x_i} (C))^d$
where
$C \rightarrow \overline{M}_{0n}$
is the universal curve,
$x_i : \overline{M}_{0n} \rightarrow C$
is the
$i$--th
section, and
$T^*_{x_i}$
is the relative cotangent sheaf to $C$ at
$x_i$.

\smallskip

In order to calculate
$\tau_1^{(i)}$,
identify
$C \rightarrow \overline{M}_{0n}$
with the morphism
$\overline{M}_{0,n+1} \rightarrow \overline{M}_{0n}$
forgetting the
$(n+1)$--th
section. Then the section
$x_i (\overline{M}_{0n})$
becomes the boundary divisor
$D_i := D_{\{i,n+1\},\{1,\ldots, \hat \iota, \ldots,n\}}$
in
$\overline{M}_{0,n+1}$,
and $\tau_1^{(i)}$ becomes the pull back of
$-D^2_i$.
Applying (1.7) to this situation we get:

\proclaim{\quad 1.6.3. Proposition}

$$\tau_1^{(i)} = \sum \Sb i\in S \subset \{1,\ldots,n\} \\ |S| \ge 2,\, n-|S|
\ge 2 \endSb \frac {(n-|S|)(n-|S|-1)} {(n-1)(n-2)} \, D_{s,\{1,\ldots,n\}
\setminus S} \ \roman{mod} \ I_S \tag 1.8$$
\endproclaim

{\bf 1.7. Multiplication formulas III.}
The functional
$\int_{\overline{M}_{0,S}} : H^*(\overline{M}_{0,S}) \rightarrow K$
is given by

$$m(\tau) \longmapsto \left\{ \aligned &1, \quad \text{if} \quad \roman{deg}\
m(\tau) =|S|-3,\\
&0 \quad \text{otherwise.} \endaligned \right.$$
Notice that
$\roman{deg} \ m(\tau) = |S|-3$
iff
$|v|=3$
for all
$v \in V_{\tau}$,
and
$\overline{M}_{\tau}$
is a point in this case. We put
$\langle \sigma_1, \, \sigma_2 \rangle =\int_{\overline{M}_{0S}} m(\sigma_1)
m(\sigma_2)$
and set to calculate this intersection index for the case when
$\roman{deg}\ m(\sigma_1) +\roman{deg}\, m(\sigma_2) =|S|-3$.
Generally, we will write
$\langle m \rangle$
instead of
$\int_{\overline{M}_{0S}}m$.
The following notions and results are due to Ralph Kaufmann. We can assume that
all pairs of different divisors of
$m(\sigma_1)$
and
$m(\sigma_2)$
are compatible, otherwise
$\langle \sigma_1,\sigma_2 \rangle =0$.
Put
$\tau =\sigma_1 \times \sigma_2$
in the category of $S$--morphisms. This is a tree with a marked subset of edges
$E$ corresponding to
$D_{\sigma}$'s
whose squares divide
$m(\sigma_1) m(\sigma_2)$.
We denote by
$\delta$
the subgraph of
$\tau$
consisting of $E$ and its vertices.

\smallskip

Consider an orientation of all edges of
$\delta$.
Call it {\it good} if for all vertices $v$ of
$\tau$,
the number of ingoing edges equals
$|v|-3$,
where
$|v|$
means the valency in
$\tau$.
Notice that for
$v \notin V_{\delta}$
we interpret this as
$|v| =3$.

\proclaim{\quad 1.7.1. Proposition (R. Kaufmann)}
There cannot exist more than one good orientation of
$\delta$.
If there is none, we have
$\langle \sigma_1,\sigma_2 \rangle =0$.
If there is one, we have

$$\langle \sigma_1,\sigma_2\rangle = \prod_{v\in V_{\tau}} (-1)^{|v|-3}
(|v|-3)! \tag 1.9$$
\endproclaim

For a proof, see Appendix.

\medskip

{\bf Remark.}
Notice that (1.9) depends only on
$\tau$
whereas $E$ influences only the existence of the good orientation. Curiously,
(1.9) coincides with the virtual Euler characteristics of the {\it non-compact}
moduli space
$M_{\tau} = \prod_{v\in V_{\tau}} M_{0,|v|}$.
We do not know why this is so.

\newpage

\centerline{\bf \S \, 2. Boundary strata and the additive structure of
$H^*(\overline{M}_{0S})$ }

\bigskip

{\bf 2.1. Basic linear relations.}
Let
$|S| \ge 4$.
Consider a system
$(\tau,v,\bar i, \bar j,\bar k, \bar l)$
where
$\tau$
is an $S$-tree,
$v \in V_{\tau}$
is a vertex with
$|v| \ge 4$
and
$\bar i, \bar j, \bar k, \bar l \in F_{\tau} (v)$
are pairwise distinct flags (taken in this order). Put
$T =F_{\tau} (v) \setminus \{ \bar i,\bar j,\bar k,\bar l\}$.
For any ordered 2-partition of $T$,
$\alpha =\{T_1,T_2\}$,
(one or both
$T_i$
can be empty) we can define two trees
$\tau^{\prime} (\alpha)$
and
$\tau^{\prime\prime}(\alpha)$.
The first one is obtained by inserting a new edge $e$ at
$v\in V$
with branches
$\{\bar i,\bar j,T_1\}$
and
$\{\bar k,\bar l,T_2\}$
at its edges. The second one corresponds similarly to
$\{\bar k,\bar j,T_1\}$
and
$\{\bar i,\bar l,T_2\}$.
We remind that
$S(\bar i)$
is the set of labels of tails belonging to the branch of
$\bar i$.

\smallskip

\proclaim {2.1.1. Proposition}
We have

$$R(\tau,v,\bar i,\bar j,\bar k,\bar l) := \sum_{\alpha}
[m(\tau^{\prime}(\alpha)) - m(\tau^{\prime\prime} (\alpha))] \equiv 0 \
\roman{mod} \ I_S \tag 2.1$$
\endproclaim

{\bf Proof.}
Choose
$i \in S(\bar i),\, j\in S(\bar j),\, k\in S(\bar k),\, l\in S(\bar l)$,
and calculate
$R_{ijkl} m(\tau) \equiv 0 \ \roman{mod} \ I_S$,
where
$R_{ijkl}$
is defined by (0.1). Consider e.g. the summands
$D_{\sigma} m(\tau)$
for
$ij\sigma kl$.

\midspace{7cm} \caption{Figure 5.}

 From the picture of
$\tau$
it is clear that
$D_{\sigma}$
does not divide
$m(\tau)$.
If
$D_{\sigma} m(\tau)$
does not vanish modulo
$I_S$,
we must have
$D_{\sigma} m(\tau) =m(\sigma \times \tau)$,
and
$\sigma \times \tau$
is of the type
$\tau^{\prime} (\alpha)$.
Similarly, the summands of
$D_{\sigma} m(\tau)$ with $kj\sigma il$ are of the type $m(\tau^{\prime\prime}
(\alpha )).$

\proclaim{\quad 2.2. Theorem}
All linear relations modulo
$I_S$
between good monomials of degree
$r+1$
are spanned by the relations (2.1) for
$|E_{\tau}|=r$.
\endproclaim

{\bf Proof.}
For
$r=0$
this holds by definition of
$I_S$.
Generally, denote by
$H_{*S}$
the linear space, generated by the symbols
$\mu (\tau)$
for all classes
$\tau$
of stable $S$-trees modulo isomorphisms, satisfying the analog of the relations
(2.1)

$$r(\tau,v,\bar i,\bar j,\bar k,\bar l) := \sum_{\alpha}
[\mu(\tau^{\prime}(\alpha)) -\mu(\tau^{\prime\prime}(\alpha))] =0 \tag 2.2$$
Denote by 1 the symbol
$\mu (\rho)$
where
$\rho$
is one-vertex tree.

\proclaim{\quad 2.2.1. Main Lemma}
There exists on
$H_{*S}$
a structure of
$H^*_S$--module given by the following multiplication formulas reproducing
(1.1), (1.3) and (1.4):

$$D_{\sigma} \mu(\tau) =\mu(\sigma \times \tau) \tag 2.3$$
if
$D_{\sigma} m(\tau)$
is a good monomial;

$$D_{\sigma} \mu(\tau) =0, \tag 2.4$$
if there exists a divisor
$D_{\sigma^{\prime}}$
of
$m(\tau)$
such that
$a(\sigma,\sigma^{\prime}) =4$;

$$D_{\sigma} \mu(\tau) =-\sum \Sb T\subset T_1 \\|T| \ge 1 \endSb \mu (tr_{T,e}
(\tau)) -\sum \Sb T\subset T_2 \\|T| \ge 1 \endSb \mu (tr_{T,e} (T)) \tag 2.5$$
if
$D_{\sigma}$
divides
$M(\tau)$,
and $e$ corresponds to
$\sigma$.
The notation in (2.5) is the same as in (1.4).
\endproclaim

{\bf Deduction of Theorem 2.2 from the Main Lemma.}
Since the monomials
$m(\tau)$
satisfy (2.1), there exists a surjective linear map
$a : H_{*S} \rightarrow H^*_S : \mu(\tau) \mapsto m(\tau)$.
On the other hand, from (2.3) it follows that
$m(\sigma) \mu (\tau) =\mu (\sigma \times \tau)$
if
$m(\sigma) m(\tau)$
is a good monomial. Hence we have a linear map
$b : H^*_S \rightarrow H_{*S} : m(\tau) \mapsto \mu(\tau) =m(\tau)1$
inverse to $a.$
Therefore
$\roman{dim} \ H_{*S} = \roman{dim} \ H^*_S$
so that the Theorem 2.2 follows.

\smallskip

We now start proving the Main Lemma.

\smallskip

{\bf 2.2.3. (2.5) is well defined.}
The r.h.s. of (2.5) formally depends on
the choice of
$\bar i,\bar j,\bar k,\bar l$.
We first check that different choices coincide modulo (2.2). It is possible to
pass from one choice to another by replacing one flag at a time. So let us
consider
$\bar i^{\prime} \ne \bar i,\bar j,\bar k,\bar l$
and write the difference of the right hand sides of the relations (2.5) written
for
$(\tau,v,\bar i,\bar j, \bar k,\bar l)$
and
$(\tau, v,\bar i^{\prime},\bar j,\bar k, \bar l)$.
The terms corresponding to those $T$ that do not contain
$\{\bar i,\bar i^{\prime}\}$
cancel. This includes all terms with
$T \subset T_2$.
The remaining sum can be rewritten as

$$-\sum_{T\subset T_1\setminus\{\bar i,\bar i^{\prime},\bar j\}} \left[
\mu(tr_{T\cup\{\bar i^{\prime}\}} (\tau)) -\mu (tr_{T\cup \{i\}} (\tau))
\right] \eqno{(2.6)}$$
where now $T$ can be empty.

\smallskip

We contend that (2.6) is of the type (2.2). More precisely, consider any of the
trees
$tr_{T\cup \{\bar i^{\prime}\}} (\tau),\, tr_{T\cup\{\bar i\} } (\tau)$
and contract the edge whose vertices are incident to the flags
$\bar i,\bar j,\bar i^{\prime}$.
We will get a tree
$\sigma$
and its vertex
$v \in V_{\tau}$.
The pair
$(\sigma,v)$
up to a canonical isomorphism does not depend on the transplants we started
with. In
$F_{\sigma} (v)$
there are flags
$\bar i,\bar j,\bar i^{\prime}$
and one more flag whose branch contains both $k$ and $l$ and which we denote
$\bar h$.
Then (2.6) is
$-r(\sigma,v,\bar i, \bar j, \bar i^{\prime},\bar h)$
(see (2.2)). This is illustrated by the Figure 6.

\newpage

\midspace{10cm} \caption{Figure 6. The edge $\ne$ is contracted to $v$}

{\bf 2.2.4. Operators $D_{\sigma}$ on $H_{*S}$ pairwise commute.}
We have to prove the identities

$$D_{\sigma_1} (D_{\sigma_2}  \mu (\tau)) =D_{\sigma_2} (D_{\sigma_1} \mu
(\tau)). \eqno{(2.7)}$$
Consider several possibilities separately.

\roster
\item "i)" {\it There exists a divisor}
$D_{\sigma}$
{\it of}
$m(\tau)$
{\it such that}
$a(\sigma_1,\sigma) =4$,
{\it so that}
$D_{\sigma_1}\mu(\tau) =0$.
\endroster

If
$D_{\sigma_2} \mu (\tau) =0$
as well, (2.7) is true. If
$D_{\sigma_2} \mu(\tau) =\mu(\sigma_2 \times \tau)$,
then
$D_{\sigma}$
divides
$m(\sigma_2 \times \tau)$,
and (2.7) is again true. Finally, if
$D_{\sigma_2}$
divides
$m(\tau)$,
then
$\sigma_2 \ne \sigma$
(otherwise
$m(\tau)$
would not be a good monomial). Hence the transplants
$tr_{T,e}(\tau)$
entering the formula of the type (2.5) which we can use to calculate
$D_{\sigma_2} \mu (\tau)$
will all contain an edge corresponding to
$\sigma$
so that
$D_{\sigma_1} (tr_{T,e} (\tau)) =0$,
and (2.7) again holds.

\smallskip

The same argument applies to the case when
$D_{\sigma_2} \mu(\tau) =0$.

\smallskip

 From now on we may and will assume that for any divisor
$D_{\sigma}$
of
$m(\tau)$
we have
$a(\sigma,\sigma_1) \le 3,\, a(\sigma,\sigma_2) \le 3$,
and that
$\sigma_1 \ne \sigma_2$.

\roster
\item"ii)" $a(\sigma_1,\sigma_2) =4$
{\it and}
$D_{\sigma_2}$
{\it divides}
$m(\tau)$.
\endroster

Then
$D_{\sigma_1}$
does not divide
$m(\tau)$,
so that
$D_{\sigma_1} \mu (\tau) =\mu (\sigma_1 \times \tau)$,
and
$D_{\sigma_2} (D_{\sigma_1} \mu(\tau)) =0$.
On the other hand,
$D_{\sigma_2} \mu (\tau)$
is a sum of transplants to the midpoint of the edge, corresponding to
$\sigma_2$.
Each such transplant has an edge giving the 2-partition
$\sigma_2$,
so that
$D_{\sigma_1} (D_{\sigma_2} \mu (\tau)) =0$.

\smallskip

The case
$a(\sigma_1,\sigma_2) =4$
and
$D_{\sigma_1}/m(\tau)$
is treated in the same way.

\smallskip

Hence from this point on we can and will in addition assume that
$a(\sigma_1,\sigma_2) =3$.

\roster
\item"iii)" $D_{\sigma_1}$
{\it does not divide}
$m(\tau)$.
\endroster

If
$D_{\sigma_2}$
does not divide
$m(\tau)$
as well, then
$D_{\sigma_1} (D_{\sigma_2} \mu (\tau)) = D_{\sigma_1} \mu (\sigma_2 \times
\tau) =\mu(\sigma_1 \times \sigma_2 \times \tau) =D_{\sigma_2} (D_{\sigma_1}
\mu (\tau))$.
If
$D_{\sigma_2}$
divides
$m(\tau)$,
we will use a carefully chosen formulas of the type (2.5) for the calculation
of
$D_{\sigma_2} \mu(\tau)$.
Namely, let
$v_1$
be the (unique) vertex of
$\tau$
which gets replaced by an edge in
$\sigma_1 \times \tau$,
and let
$e_2$
be the edge of
$\tau$
corresponding to
$D_{\sigma_2}$.
Let
$u_2,\, u_1$
be the vertices of
$e_2$
such that
$u_1$
can be joined to
$v_1$
by a path not passing by
$e_2$.

\smallskip

Consider first the subcase
$u_1 \ne v_1$.
Choose some
$\bar i,\bar j \in F_{\tau} (u_2)$
and
$\bar k,\bar l \in F_{\tau} (u_1)$
in such a way that
$\bar l$
starts a path leading from
$u_1$
to
$v_1$.
Use these
$\bar i,\bar j,\bar k,\bar l$
in a formula of the type (2.5) to calculate
$D_{\sigma_2} \mu(\tau)$
and then
$D_{\sigma_1} (D_{\sigma_2} \mu (\tau))$,
that will insert an edge instead of the vertex
$v_1$
which survives in all the transplants entering
$D_{\sigma_2} \mu(\tau)$.
Then calculate
$D_{\sigma_2} (D_{\sigma_1}\mu(\tau))$
by first inserting the edge at
$v_1$,
and then constructing the transplants not moving
$\bar i, \bar j, \bar k, \bar l$.
Since by our choice of
$\bar l$
we never transplant the branch containing
$v_1$,
the two calculations will give the same result.

\smallskip

Now let
$v_1 =u_1$.
Let
$\{S_1,S_2\}$
be the 2-partition of $S$ corresponding to
$\sigma_1$.
Since
$\sigma \times \tau$
exists,
$\{S_1,S_2\}$
is induced by a partition of
$F_{\tau}(v_1) = \bar S_1 \coprod \bar S_2$.
We denote by
$\bar S_2$
the part to which the flag
$(v_1 =u_1,\, e_2)$
belongs. Let
$\bar T =\bar S_2 \setminus \{ (v_1 =u_1,e_2)\}$.
This set is non-empty because otherwise
$e_2$
would correspond to
$\{S_1,S_2\}$
and we would have
$\sigma_1 =\sigma_2$.
Take
$\bar i,\bar j \in F_2(\tau),\, \bar k \in \bar S_1$
and
$\bar l \in T$:
see Figure 7.

\midspace{8cm} \caption{Figure 7}

Now consider
$D_{\sigma_2} (D_{\sigma_1} \mu(\tau))$
and
$D_{\sigma_1} (D_{\sigma_2} \mu (\tau))$.
To calculate the first expression we form a sum of transplants of
$\sigma_1 \times \tau$.
To calculate the second one, we form transplants of
$\tau$,
and then insert an edge at
$v_1 =u_1$.

\smallskip

The transplants corresponding to the branches at
$u_2$
will be the same in both expressions. The transplants corresponding to the
subsets
$T \subset \bar T \setminus \{\bar l\}$
will also be the same. In addition, the second expression will contain the
terms
$-D_{\sigma_1} (\mu (tr_{T,e_2}(\tau)))$
where
$T \cap \bar S_1 \ne \emptyset$.
But each such term vanishes. In fact, consider the 2--partition
$\rho =\{R_1,R_2\}$
of $S$ corresponding to the edge of
$tr_{T,e_2} (\tau)$
containing the flag
$(v_1 =u_1,e_2)$,
and let
$k,l \in R_1$.
A glance to the third tree of the Figure 7 shows that
$a(\rho,\sigma_1)=4$,
because if
$\bar t \in T \cap \bar S_1,\, t \in S(\bar t)$,
then
$kt\sigma_1 il$
and
$kl\rho it$.
Hence the extra terms are irrelevant.

\smallskip

The case when
$D_{\sigma_2}$
does not divide
$m(\tau)$
is treated in the same way. It remains to consider the last possibility.

\smallskip

\roster
\item"iv)" $D_{\sigma_1}$
{\it and}
$D_{\sigma_2}$
{\it divide}
$m(\tau),\, a(\sigma_1,\sigma_2)=3$.
\endroster

Denote by
$e_1$
(resp.
$e_2$)
the edge corresponding to
$\sigma_1$
(resp.
$\sigma_2$).
Let
$u_1,u_2$
(resp.
$v_1,v_2)$
be the vertices of
$e_1$
(resp.
$e_2)$
numbered in such a way that there is a path from
$u_2$
to
$v_1$
not passing through
$e_1, \, e_2$
(the case
$u_2 =v_1$
is allowed). To calculate the multiplication by
$D_{\sigma_1}$
choose
$\bar i,\bar j \in F_{u_1} (\tau) \setminus \{(u_1,e_1)\},\, \bar l$
on the path from
$u_2$
to
$v_1$ if
$u_2 \ne v_1$,
and
$\bar l =(v_1,e_2)$
if
$u_2 =v_1;\, \bar k \in F_{\tau}(v_2) \setminus \{\bar l\}$.
To calculate the product by
$D_{\sigma_2}$,
choose similarly
$\bar k^{\prime},\bar l^{\prime} \in F_{\tau} (v_2)\setminus \{(v_2,e_2)\}$,
$\bar i^{\prime} \in F_{\tau} (v_1)$
on the path from
$v_1$
to
$u_2$,
if
$v_1 \ne u_2$,
and
$\bar i^{\prime} =(u_2,e_1)$
if
$v_1 =u_2,\, \bar j^{\prime} \in F_{\tau} (v_1)$
(see Figure 8).

\midspace{4cm} \caption{Figure 8}

The critical choice here is that of
$\bar l$
and
$\bar i^{\prime}$.
It ensures that calculating
$D_{\sigma_1} (D_{\sigma_2} \mu(\tau))$
and
$D_{\sigma_2} (D_{\sigma_1} \mu (\tau))$
we will get the same sum of transplanted trees. This ends the proof of (2.7).

\medskip

{\bf 2.2.5. Compatibility with $I_S$-generating relations.}
If
$D_{\sigma_1} \, D_{\sigma_2} =0$
because
$a(\sigma_1,\sigma_2) =4$,
one sees that
$D_{\sigma_1} (D_{\sigma_2} \mu (\tau)) =0$
looking through various subcases in 2.2.4. It remains to show that
$R_{ijkl} \mu (\tau) =0$
where
$R_{ijkl}$
is defined by (0.10).

\smallskip

Consider the smallest connected subgraph in
$\tau$
containing the flags
$i,j,k,l$.
The Figure 9 gives the following exhaustive list of alternatives. Paths from
$i$ to $j$ and from $k$ to $l$: i) have at least one common edge; ii) have
exactly one common vertex; iii) do not intersect.

\midspace{6cm} \caption{Figure 9}

Consider them in turn.

\smallskip

i). Let $e$ be an edge common to the paths
$ij$
and
$kl$.
Denote by
$\rho$
the respective 2-partition. Then
$ik\rho jl$
or
$il\rho kj$.
Therefore any summand of
$R_{ijkl}$
annihilates $D_{\rho}$ so that $R_{ijkl}\mu (\tau )=0$ in
view of (2.4).

\smallskip

ii). Let $v$ be the vertex common to the paths $ij$ and $kl$.
Then exactly the same calculation as in the proof
of the Proposition 2.1.1 shows that
$$
R_{ijkl}\mu(\tau )=\sum_{\alpha}
[\mu (\tau^{\prime}(\alpha)) -\mu (\tau^{\prime\prime}(\alpha))]=0
$$
(notation as in (2.1) and (2.2)).

\smallskip

iii). This is the most complex case. Let us draw a more detailed
picture of $\tau$ in the neighborhood of the subgraph we are
considering (Figure 10).

\midspace{6cm} \caption{Figure 10}

Let $v_1$ be the vertex on the path $ij$ which is connected
by a sequence of edges $e_1,\dots ,e_m$ ($m\ge 1$) with the vertex
$v_m$ on the path $kl$ so that $e_a$ has vertices $(v_a,v_{a+1})$
in this order. Let $T_a$ be the set of flags at $v_a$ which do
not coincide with $\bar i,\bar j,\bar k,\bar l,$ and do not belong
to $e_{a-1},e_{a}.$

\smallskip

Consider any summand $D_{\sigma}$ of $R_{ijkl}.$ If
$jk\sigma il,$ then $D_{\sigma}\mu (\tau )=0$ because each
edge $e_a$ determines a partition $\rho$ of $S$ such that
$ij\rho kl.$ From now on we assume that $ij\sigma kl.$
Then $D_{\sigma}\mu (\tau )$ can be nonzero if one of the
two alternatives holds:

\smallskip

a). For some $v_a$, there exists a partition
$T_a=T_a^{\prime}\coprod T_a^{\prime\prime},$(with $|T_a^{\prime}|\ge 1,
|T_a^{\prime\prime}\ge 1,$ except for the case $a=1$ where $T_1^{\prime}$
can be empty, and $a=m$ where $T_m^{\prime}$ can be empty) such that the
following two sets
$$
S_1=S(\bar i)\coprod S(\bar j)\coprod S(T_1^{\prime})\coprod \dots
\coprod S(T_a^{\prime}),
$$
$$
S_2=S(T_a^{\prime\prime})\coprod S(T_{a+1})\coprod \dots
\coprod S(T_m)\coprod S(\bar k)\coprod S(\bar l)
$$
form the 2-partition corresponding to $\sigma .$ In this case
$$
D_{\sigma}\mu (\tau )=\mu (\sigma\times\tau ),
$$
and $\sigma\times\tau$ is obtained by inserting a new edge at $v_a$
and by distributing $T_a^{\prime}$ and $T_a^{\prime\prime}$
at different vertices of this edge.

\smallskip

b). For some $e_a,$ the two sets
$$
S_1=S(\bar i)\coprod S(\bar j)\coprod (\coprod_{i\le a}S(T_i)),
$$
$$
S_2=(\coprod_{i\ge a+1}S(T_i))\coprod S(\bar k)\coprod S(\bar l)
$$
form the 2-partition corresponding to $\sigma .$

\smallskip

In this case $D_{\sigma}$ divides $m(\tau ),$ and in order
to calculate $D_{\sigma}\mu (\tau )$ using a formula of the type
(2.5) we must first choose two pairs of flags at two vertices
of $v_a .$

\smallskip

Contributions from a) and b) come with opposite signs, and we contend that they
completely cancel each other.

\smallskip

To see the pattern of the cancellation look first at the case a) at $v_1.$
It brings (with positive sign) the contributions corresponding
to the following trees. Form all the partitions $T_1=T_1^{\prime}\coprod
T_1^{\prime\prime}$ such that $T_1^{\prime\prime}\ne \emptyset ,$
where $T_1=F_{\tau}(v_1)\setminus \{ \bar i,\bar j,(v_1,e_1)\}.$
Transplant all $T_1^{\prime\prime}$--branches to the midpoint of $e_1.$
Denote the new vertex $v_1^{\prime}.$ The result is drawn
as Figure 11.

\midspace{6cm} \caption{Figure 11}

Now consider the terms of the type b) for the edge $e_1.$
If $m=2,$ we choose for the calculation of $D_{\sigma_1}\mu (\tau )$
(where $\sigma_1$ corresponds to $e_1$) the flags $\bar i,\bar j,\bar k,\bar
l.$
If $m> 2,$ we choose the flags $\bar i,\bar j,(v_2,e_1),t\in T_2.$
Then we get the sum of two contributions. One will consist
of the trees obtained by transplanting branches at $v_1.$
They come with negative signs and exactly cancel the previously considered
terms of the type a). If $m=2,$ the second group will
cancel the terms of the type a) coming from $v_2.$

\smallskip

Consider a somewhat more difficult case $m>2.$ Then this second group
of terms comes from the trees indexed by the partitions
$T_2=T_2^{\prime}\coprod T_2^{\prime\prime}, t\in T_2^{\prime\prime},
T_2^{\prime}\ne \emptyset .$ Branches corresponding to
$T_2^{\prime}$ are transplanted to the midpoint $v_1^{\prime}$ of the
edge $e_1.$ These terms come with negative signs: see Figure 12.

\midspace{6cm} \caption{Figure 12}

These trees in turn cancel with those coming from the terms
of the type a) at the vertex $v_2$ with positive sign.
However, there will be additional terms of the type a) for which
$t\in T^{\prime}_2.$ They will cancel with one group of
transplants contributing to
$D_{\sigma_2}\mu (\tau )$ where $\sigma_2$ corresponds to the edge
$e_2$ of the Figure 10, if for the calculation of $D_{\sigma_2}\mu (\tau )$
one uses (2.5) with the following choice of flags:
$(v_2,e_1),t$ at one end, $(v_3,e_1)$, some $t^{\prime}\in T_3$
at the other end (this last choice must be replaced by
$\bar k,\bar l,$ if $m=3$).

\smallskip

The same pattern continues until all the terms cancel.

\medskip

{\bf 2.2.6. Compatibility with relations (2.2).} By this time
we have checked that the action of any element of $F_S/I_S$
on the individual generators $\mu (\tau )$ of $H_{*S}$ is well defined
modulo the span $I_{*S}$ of relations (2.2). It remains to show
that the subspace in $\oplus_{\tau}K\mu (\tau )$ spanned by these relations
is stable with respect to this action. But the calculation in the
proof of the Proposition 2.2.1 shows that
$$
r(\tau,v,\bar i,\bar j,\bar k,\bar l) \equiv
m(\tau )r_{ijkl}\ \roman{mod}\ I_{*S},
$$
where $r_{ijkl}$ is obtained from $R_{ijkl}$ by replacing
$m(\sigma )$ with $\mu (\sigma )$. To multiply this by any element
of $H^*_S$ we can first multiply it by $m(\tau )$, then
represent the result as a linear combination of good monomials,
and finally multiply each good monomial by $r_{ijkl}.$ The result
will lie in $I_{*S}.$

\smallskip

This finishes the proof of the Main Lemma and the Theorem 2.2.

\newpage

\centerline{\bf \S 3. Cohomological Field Theories of rank 1.}

\bigskip

{\bf 3.1. Notation.} Rank of a CohFT on $(H,g)$ is the
(super)dimension of $H.$ In this section we consider the case
$\roman{dim}\ H=1.$ To slightly simplify notation let us assume that
all square roots exist in $K.$ Then $H=K\Delta_0,\ g(\Delta_0,\Delta_0)=1,\
\Delta = \Delta_0^{\otimes 2}\in H^{\otimes 2}.$ The basic vector
$\Delta_0$ is defined up to a sign. We will consider its choice as
a rigidification and without further ado call such
rigidified theories simply CohFT's of rank one.

\smallskip

A structure of CohFT on $(H,\Delta_0)$ boils down to
a sequence of cohomology classes (generally non--homogeneous)
$$
c_n:=I_n(\Delta_0^{\otimes n})\in H^*(\overline{M}_{0n},K)^{S_n},
\ n\ge 3,
\eqno{(3.1)}
$$
satisfying the identities
$$
\varphi^*_{\sigma}(c_n)=c_{n_1+1}\otimes c_{n_2+1},
\eqno{(3.2)}
$$
where $\phi_{\sigma}:\ \overline{M}_{0,n_1+1}\times\overline{M}_{0,n_2+1}
\to \overline{M}_{0n}$ is the embedding of the boundary divisor
corresponding to a partition $\sigma$ (see (0.3)). Put
$c_n=\sum_{i=0}^{n-3}c^{(i)}_n,\ c^{(i)}_n\in H^{2i}(\overline{M}_{0n}).$
Changing sign of $\Delta_0$ leads to $c_n\mapsto (-1)c_n.$

\smallskip

The tensor product formula (0.8) becomes
$$
\{c^{\prime}_n\}\otimes\{c^{\prime\prime}_n\}=
\{c^{\prime}_n\wedge c^{\prime\prime}_n\},
\eqno({3.3)}
$$
if we agree that $(H^{\prime},\Delta_0^{\prime})\otimes
(H^{\prime\prime},\Delta_0^{\prime\prime})=
(H^{\prime}\otimes H^{\prime\prime},\Delta_0^{\prime}
\otimes\Delta_0^{\prime\prime}).$

\smallskip

Here are some simple consequences of (3.2) and (3.3).

\smallskip

\proclaim{\quad 3.1.1} The theory $c_n=c_n^{0}=[\overline{M}_{0n}]$
for $n\ge 3$ is the identity with respect to the tensor product.
\endproclaim

\smallskip

\proclaim{\quad 3.1.2} The theories $c_n(t)=t^{n-2}[\overline{M}_{0n}],$
$t\in K^*$, form a group isomorphic to $K^*.$
\endproclaim

\smallskip

\proclaim{\quad 3.1.3} The theory $\{c_n\}$ is invertible iff
$c^{(0)}_3\ne 0.$ Any invertible theory is a tensor
product of one of the type $\{c_n(t)\}$ and one with
$c_n^{(0)}=1$
 for all $n$, and this decomposition
is unique.
\endproclaim

\smallskip

\proclaim{\quad 3.1.4} Assume that $c_n^{(0)}=1$ and put
$\lambda_n=\roman{log}\ c_n\in H^*(\overline{M}_{0n},K)^{S_n}.$
Then (3.2) becomes
$$
\varphi^*_{\sigma}(\lambda_n)=\lambda_{n_1+1}\otimes 1+
1\otimes \lambda_{n_2+1},
\eqno{(3.4)}
$$
and (3.3) becomes
$$
\{\lambda_n^{\prime}\}\otimes\{\lambda_n^{\prime\prime}\}=
\{\lambda_n^{\prime}+\lambda_n^{\prime\prime}\}.
\eqno{(3.5)}
$$
Vice versa, any sequence of classes
$\lambda_n\in H^*(\overline{M}_{0n},K)^{S_n}$
satisfying (3.4) gives rise to a CohFT of rank 1,
$c_n=\roman{exp}\,\lambda_n.$ We can say that
$\{\lambda_n\}$ forms a logarithmic CohFT of rank 1.
\endproclaim

\smallskip

\proclaim{\quad 3.1.5} There is a canonical bijection between
the set of the isomorphism classes of CohFT's of rank 1 and the set
of infinite sequences $(C_3,C_4,\dots )\in K^{\infty}$ given by
$$
C_n=\int_{\overline{M}_{0n}}c_n=\int_{\overline{M}_{0n}}c_n^{(n-3)}.
$$
\endproclaim

\smallskip

In fact, this is a particular case of the equivalence stated in
0.4 and 0.5, because any formal series in one variable
$\Phi (x)=\sum \frac{C_n}{n!}x^n$ satisfies the associativity
equations.

\smallskip

Formula (0.7) reconstructing $c_n^{(i)}$ from $C_m$ becomes
$$
\forall \,\tau, |E_{\tau}|=i:\quad
\int_{\overline{M}_{0n}}[m(\tau )]\wedge c_n^{(i)}=
\prod_{v\in V_{\tau}}C_{|v|},
\eqno{(3.6)}
$$
where $[m(\tau )]\in H^*(\overline{M}_{0n})$ is the image of
the good monomial $m(\tau ).$ We do not know
nice formulas for the tensor product in terms of
coordinates $(C_i)$. The main goal of this section
is to show that there are natural coordinates defined
geometrically that are simply additive with respect
to the tensor multiplication of invertible theories.
This is a reformulation of certain identities from [AC].

\medskip

{\bf 3.2. Mumford classes.} Consider the universal curve $p_n:\ X_n\to
\overline{M}_{0n}$ and its structure sections $s_i:\ \overline{M}_{0n}
\to X_n,\ i=1,\dots ,n.$ Let $x_i\subset X_n$ be the image of $s_i$,
$\omega$ the relative dualizing sheaf on $X_n.$ For $a=1,2,\dots $
put
$$
\omega_n(a):=p_{n*}\left( c_1(\omega (\sum_{i=1}^n x_i))^{a+1}\right)
\in H^{2a}(\overline{M}_{0n},\bold{Q})^{S_n}.
\eqno{(3.9)}
$$

\smallskip

It is proved in [AC] that for any $a\ge 1$ (in fact, $a=0$ as well)
$\{\omega_n(a)\,|\,n\ge 3\}$ satisfy (3.4) i.e., form a logarithmic
field theory. Hence we can construct an infinite--dimensional
family of invertible theories of rank one:
$$
\omega_n[s_1,s_2,\dots ]:=\roman{exp}\
(\sum_{a=1}^{\infty}s_a\omega_n(a)),\ n\ge 3.
\eqno{(3.8)}
$$

\medskip

\proclaim{\quad 3.2.1. Theorem} $(s_a)$ form a coordinate system on
the space of isomorphism classes of theories with $c^{(0)}_3=1$,
defining its group isomorphism with $K^{\infty}_+$.
\endproclaim

{\bf Proof.} The sum in the r.h.s. effectively stops at $a=n-3$
(cf. (3.7)). The $C_n$--coordinate of the theory (3.8) is therefore
$$
\int_{\overline{M}_{0n}}\omega_n[s_1,s_2,\dots ]=
s_{n-3}\int_{\overline{M}_{0n}}\omega_n(n-3)+
P_n(s_1,\dots ,s_{n-4}),\ n\ge 4,
$$
where $P_n$ is a universal polynomial. Hence it remains to check that the
coefficient at $s_{n-3}$ does not vanish. But this follows from the well known
fact that $\omega(\sum_{i=1}^nx_i)$
is an ample sheaf on $X_n.$

\medskip

{\bf 3.2.2. Remark.} The theories we are considering here are tree
level ones in the terminology of [KM]. The general definition
of a CohFT given there involves maps $I_{g,n}:\ H^{\otimes n}\to
H^*(\overline{M}_{g,n})$ for all stable pairs $(g,n).$
The classes $\omega_n(a)$ given by (3.7) can be automatically
defined in this larger generality, and the extension of
the property (3.4) is proved in [AC] for all $(g,n).$
Therefore formulas (3.8) in fact define full (any genus)
rank one theories.

\smallskip

However we do not know whether Theorem 3.2.3 extends to the general case
because
it is unclear whether functions $C_{g,n}:=\int_{\overline{M}_{g,n}}c_{g,n}$
form a coordinate system on the space of full rank 1 theories.
In fact, they probably do not, because of the presence of
non--trivial cusp classes in $H^*(\overline{M}_{g,n})$ having vanishing
restrictions to the boundary.

\medskip

{\bf 3.3. Potential of rank 1 theories.} Denote by $\Phi (x;s_1,s_2,\dots )$
the potential of $\omega_n[s_1,s_2,\dots]$ at $x\Delta$ (see (0.6)).
We have from (3.8):
$$
\Phi (x;s_1,s_2,\dots )=\sum_{n=3}^{\infty}\frac{x^n}{n!}
\sum_{(m_a):\ \sum am_a=n-3}\int_{\overline{M}_{0n}}
\prod_a \omega_n(a)^{m_a}\prod_a\frac{s_a^{m_a}}{m_a!}.
\eqno{(3.9)}
$$
We expect that (3.9) satisfies some interesting differential equations
encoding recursive relations between the numbers
$$
\int_{\overline{M}_{0n}}\prod_a \omega_n(a)^{m_a},\
\sum am_a=n-3.
$$
Some partial results are given below.

\smallskip

It is even possible that such equations for arbitrary genus are
implicit in the relations (conjectured by Witten and proved by
Di Francesco, Itzykson, and Zuber) between the numbers
denoted in [AC]
$$
\int_{W_{(m_a),n}}\prod_a\psi_i^{n_a},
$$
where $W_{(m_a),n}$ are certain combinatorial classes defined in terms of
ribbon graphs. In fact, it is conjectured in [AC] that
the dual cohomology classes of $W_{(m_a),n}$ can be expressed as
$$
\prod_{a\ge 2} \omega_{n}(a-1)^{m_a}\frac{(2^a(2a-1)!!)^{m_a}}{m_a!}
$$
plus terms of lower order and boundary terms (our $\omega_n(a)$
are denoted $k_a$ in [AC]).

\smallskip

We hope to return to this problem elsewhere. Here we will treat
the case when only one of the coordinates $s_a$ is non--zero.

\medskip

{\bf 3.4. Weil--Petersson theory.} The noncompact moduli spaces $M_{0n}$
possess a canonical Weil--Petersson hermitean metric. It is
singular on the boundary, but its K\"ahler form extends to a
closed $L^2$--current on $\overline{M}_{0n}$ thus defining a real
cohomology class $\omega^{WP}_n\in H^2(\overline{M}_{0n})^{S_n}$
(see [W] and [Z]).

\smallskip

In [AC] this class is identified as
$$
\omega^{WP}_n=2\pi^2\omega_n(1).
\eqno{(3.10)}
$$
The additivity property (3.4) for $\omega_n^{WP}$ was used by P. Zograf
([Z]) in order to calculate the WP--volumes of $\overline{M}_{0n}$.
In our framework, he calculated the coefficients of the potential
of the theory $\{\roman{exp}(\pi^{-2}\omega_n^{WP})\} .$

\medskip

{\bf 3.5. Weil--Petersson potential.} The $C_n$--coordinate of $\{\roman{exp}
(\pi^{-2}\omega_n^{WP})\}$ is
$$
\frac{1}{\pi^{2(n-3)}}\int_{\overline{M}_{0n}}
\frac{(\omega_n^{WP})^{n-3}}{(n-3)!}=
\frac{v_n}{(n-3)!}
$$
in the notation of [Z]. P. Zograf proved that $v_4=1,v_5=5,
v_6=61, v_7=1379,$ and generally
$$
v_n=\sum_{i=1}^{n-3}\frac{i(n-i-2)}{n-1}
{{n-4}\choose{i-1}}{{n}\choose{i+1}}v_{i+2}v_{n-i},\ n\ge 4.
\eqno{(3.11)}
$$
The potential (0.6) of the theory is therefore
$$
\Phi^{WP}(x):=\sum_{n=3}^{\infty}\frac{v_n}{n!(n-3)!}x^n.
$$
We can rewrite (3.11) as a differential equation for
$\Phi^{WP}(x).$ Following [M], put
$$
g(x)=x^2\frac{d}{dx}(x^{-1}\frac{d}{dx}\Phi^{WP}(x)).
$$
Then we have
$$
x(x-g)g^{\prime\prime}= x(g^{\prime})^2+(x-g)g^{\prime}.
$$

\medskip

{\bf 3.6. A generalization of Zograf's recursive relations.} Put
$$
z_n:=\int_{\overline{M}_{0n}}\omega_n(n-3),\qquad n\ge 3.
\eqno{(3.12)}
$$
Define for each $n\ge 3$ an $S_n$--invariant function $A_n$ on the
set of isomorphism classes of $n$--trees $\sigma$ with the following
property:
$$
\langle\sum_{\sigma :|E_{\sigma}|=a}A_n(\sigma )\sigma,\tau\rangle=
\cases 1\ \roman{if}\  |E_{\tau}|=n-3-a\ \roman{and}\
\exists v\in V_{\tau},\  |v|=a+3,\\
0\ \roman{otherwise,}\endcases \eqno{(3.13)}
$$
where $\langle\sigma ,\tau \rangle$ is defined by (1.9).
Presumably, (3.12) can be calculated inductively using the definition (3.7).
The solvability of (3.13) will be shown below.
To find an explicit solution one has to invert the Poincar\'e
pairing matrix restricted to the $S_n$--invariant part
of $H^*(\overline{M}_{0n}).$ Hopefully, a version of the Proposition 1.7.1
can be used to do this.

\smallskip

Put
$$
\Omega_n(a)=\cases \int_{\overline{M}_{0n}}\omega_n(a)^{\frac{n-3}{a}},\
\roman{if}\  a/(n-3),\\
0 \ \roman{otherwise.}\endcases \eqno{(3.14)}
$$
This is a part of the coefficients in (3.9).

\medskip

\proclaim{\quad 3.6.1. Theorem} For a fixed $a\ge 1$, the sequence
$\{\Omega_n(a)\},\ n\ge 3$ satisfies the recursive relations
$$
\Omega_n(a)=z_{a+3}\sum_{{n-\roman{trees}\ \sigma :}\atop{|E_{\sigma}|=a}}
A_n(\sigma)
\frac{\left(\frac{n-3}{a}-1\right) !}{\prod_{v\in V_{\sigma}}
\left(\frac{|v|-3}{a}  \right) !}
\prod_{v\in V_{\sigma}}\Omega_{|v|}(a). \eqno{(3.15)}
$$

\endproclaim

{\bf Proof.} First of all, we can calculate $\omega_n(a)$
as a functional on the homology classes of the strata
$\overline{M}_{\tau}$ where $\tau$ runs over stable
$n$--trees with $|E_{\tau}|=n-3-a.$ Namely, since
$\{\omega_n(a)\}$ form a logarithmic field theory, we have
$$
\int_{\overline{M}_{\tau}}\omega_n(a):=\langle
\varphi_{\tau}^*(\omega_n(a))\rangle=
\langle\sum_{v\in V_{\tau}}\roman{pr}_v^*(\omega_{|v|}(a))\rangle ,
\eqno{(3.16)}
$$
where $\roman{pr}_v^*:\  \overline{M}_{\tau}\to \overline{M}_{0,F_{\tau}(v)}$
is the canonical projection.

\smallskip

In (3.16), only the summands with $|v|=a+3$ can be non--vanishing,
and there can exist at most one such summand, because
$|V_{\tau}|=|E_{\tau}|+1=n-2-a,$ so that
$$
\sum_{v\in V_{\tau}}(|v|-3)=|T_{\tau}|+2|E_{\tau}|-3|V_{\tau}|=a.
$$
It follows that $\omega_n(a)$ is dual to the class
$z_{a+3}\sum_{\sigma :|E_{\sigma}|=a}A_n(\sigma )[\overline{M}_{\sigma}],$
in the notation of (3.13).

\smallskip

Similarly, one can calculate $\omega_n(a)^{\frac{n-3}{a}-1}$ as
a functional on the classes of the tree strata $[\overline{M}_{\sigma}]$
with $|E_{\sigma}|=n-3-a(\frac{n-3}{a}-1)=a$. We have, putting
$m=\frac{n-3}{a}-1:$
$$
\int_{\overline{M}_{\sigma}}\varphi_{\sigma}^*(\omega_n(a)^m)=
\langle\left(\sum_{v\in V_{\sigma}}
\roman{pr}_v^*(\omega_{|v|}(a))\right)^m\rangle=
$$
$$
\sum_{{(m_v|v\in V_{\sigma}):}\atop{\sum m_v=m}}
\frac{m!}{\prod m_v!}\prod_{v\in V_{\sigma}}
\int_{\overline{M}_{0,F_v(\sigma )}}\left(\omega_{|v|}(a))\right)^{m_v}=
$$
$$
\frac{m!}{\prod_v \left(\frac{|v|-3}{a} \right)!}
\prod_{v\in V_{\sigma}} \Omega_{|v|}(a),
$$
because only one summand, with $m_v=\frac{|v|-3}{a}$ for all $v$,
can be non--vanishing.

\smallskip

In view of (3.13), this is equivalent to (3.15), because
$$
\Omega_n(a)= \int_{\overline{M}_{0n}}\omega_n(a)^{\frac{n-3}{a}}
=\int_{\overline{M}_{0n}}\omega_n(a)^{\frac{n-3}{a}-1}\wedge
\omega_n(a),
$$
which is $\omega_n(a)^{\frac{n-3}{a}-1}$ integrated along
$$
z_{a+3}\sum_{\sigma :|E_{\sigma}|=a}A_n(\sigma )[\overline{M}_{\sigma}].
$$

{\bf 3.6.2. Remark.} Zograf's argument essentially coincides
with our reasoning for the case $a=1$.

\smallskip

It can be directly generalized to obtain more general recursive relations
for all coefficients in (3.9). However, their usefulness
depends on the understanding of $\{A_n(\sigma )\}$.

\medskip

{\bf 3.7. Twisting.} For any CohFT on $(H,g)$, we can define
a new theory tensor multiplying it by $\omega_.[s_1,s_2,\dots ].$
It would be interesting to study the dependence of its potential
on $s_1,s_2,\dots$ . This could clarify the analytic properties
of the initial theory.

\smallskip

If the initial theory corresponds to a system of GW-classes, as in [KM],
it satisfies a number of additional axioms. In particular,
it has a scaling group related to the grading of $H$, and an identity
in the quantum cohomology ring. Twisting generally destroys
these additional structures.

\newpage

\centerline{\bf Appendix: Proof of the Proposition 1.7.1.}

\smallskip

\centerline{R. Kaufmann}

\bigskip

We keep notation of sec. 1.7.

\smallskip

Consider the canonical embedding $\varphi_{\tau}:\
\overline{M}_{\tau}\to \overline{M}_{0S}.$ We start with the formula
$$
\langle \sigma_1,\sigma_2\rangle=
\langle\prod_{e\in E}\varphi_{\tau}^*(D_{\sigma (e)})\rangle ,
\eqno{(A.1)}
$$
where the cup product in the r.h.s. is taken in
$H^*(\overline{M}_{\tau})\cong \otimes_{v\in V_{\tau}}
H^*(\overline{M}_{0,F_{\tau}(v)}).$ Applying an appropriate version
of the formulas (1.4) and (1.6) we can write for any
$e\in E$ with vertices $v_1,v_2$:
$$
\varphi_{\tau}^*(D_{\sigma (e)})=-\Sigma_{v_1,e} -\Sigma_{v_2,e},
\eqno{(A.2)}
$$
where
$$
\Sigma_{v_i,e}\in H^*(\overline{M}_{0,F_{\tau}(v_i)})\otimes
\prod_{v\ne v_i}[\overline{M}_{0,F_{\tau}(v)}]
\eqno{(A.3)}
$$
and $[\overline{M}_{0,F_{\tau}(v)}]$ is the fundamental class.
Later we will choose an expression for $\Sigma_{v_i,e}$ depending
on the choice of flags denoted $\overline{i},\overline{j}$ or
$\overline{k},\overline{l}$ in (1.4).

\smallskip

Inserting (A.2) into (A.1), we get
$$
\langle \sigma_1,\sigma_2\rangle=\sum_h
\langle\prod_{{(v,e)\in F_{\delta}}\atop{h(e)=v}}(-\Sigma_{v,e})\rangle ,
\eqno{(A.4)}
$$
where $h$ runs over all orientations of $E$ considered as a choice,
for every $e\in E$, of a vertex $h(e)$ of $e.$

\smallskip

The summand of (A.4) corresponding to a given $h$ can be non--zero
only if for every $v\in V_{\delta}$ the number of factors
$(v,e)$ with $h(e)=v$ equals $\roman{dim}\ \overline{M}_{0,F_{\tau}(v)}=
|v|-3.$ This is what was called a good orientation.

\smallskip

Assume that there are two good orientations $h,h^{\prime}$ of $\delta .$
Consider the union of all closed edges on which $h\ne h^{\prime}.$
Each connected component of this union is a tree. Choose an end edge
$e$ of this tree and an end vertex $v$ of $e.$ At $v,$ the number
of $h$--incoming and $h^{\prime}$--incoming edges must coincide, but on $e$
these orientations differ. Hence there must exist an edge $e^{\prime}\ne e$
incident to $v$ upon which $h$ and $h^{\prime}$ differ. But this
contradicts to the choice of $v$ and $e$.

\smallskip

Now assume that one good orientation $h$ exists. We can rewrite (A.4) as
$$
\langle \sigma_1,\sigma_2\rangle=
\langle\prod_{v\in V_{\tau}}
\prod_{e:{h(e)=v}}(-\Sigma_{v,e})\rangle .
\eqno{(A.5)}
$$
In view of (A.3), this expression splits into a product of
terms computed in all $H^*(\overline{M}_{0,F_{\tau}(v)}), v\in V_{\tau}$
separately. Each such term depends only on $|v|,$ and we want
to demostrate that it equals $(-1)^{|v|-3}(|v|-3)!$. Put $|v|=m+3$.
We may and will assume that $m\ge 1,$ the case $m=0$ being
trivial.

\smallskip

Let us identify $F_{\tau}$ with $\{ 1,\dots ,m+3\}$ in such a way
that flags $1,\dots ,m$ belong to the edges $e$ with $h(e)=v.$
Denote by $D_{\rho}^{(m+3)}$ the class of a boundary divisor
in $H^*(\overline{M}_{0,m+3})$ corresponding to a stable
partition $\rho$ of $\{1,\dots ,m+3\}.$ We will choose flags
$m+1, m+2$ to play the role of $\overline{i},\overline{j}$
in (1.4) for any $i\in \{1,\dots ,m\}$ corresponding to an edge
$e$ in (A.5) ($v$ being now fixed), so that
the contribution of $v$ in (A.5) becomes
$$
\prod_{i=1}^m \left(
\sum_{\rho:\ i\rho\{m+1,m+2\}}-D^{(m+3)}_{\rho}\right):=g(m).
\eqno{(A.6)}
$$
We will calculate (A.6) inductively. Consider the projection map
(forgetting the $(m+3)$--th point) $p:\ \overline{M}_{0,m+3}\to
\overline{M}_{0,m+2}$ and the $i$--th section map
$x_i:\ \overline{M}_{0,m+2}\to\overline{M}_{0,m+3}$
obtained via the identification of $\overline{M}_{0,m+3}$ with the
universal curve. We have $p\circ x_i=\roman{id},$ and
$x_i$ identifies $\overline{M}_{0,m+2}$ with $D^{(m+3)}_{\sigma_i}$
where
$$
\sigma_i=\{\{i,m+3\},\{1,\dots ,\widehat{i},\dots ,m+2\}\}.
$$

Therefore
$$
\sum_{\rho:\ i\rho\{m+1,m+2\}}-D^{(m+3)}_{\rho}=
-p^*\left(\sum_{\rho^{\prime}:\ i\rho^{\prime}\{m+1,m+2\}}
D^{(m+2)}_{\rho^{\prime}}\right)
-x_{i*}([\overline{M}_{0,m+2}]),
\eqno{(A.7)}
$$
where $\rho^{\prime}$ runs over stable partitions of
$\{1,\dots ,m+2\}.$ We now insert (A.7) into (A.6) and
represent the resulting expression as a sum of products
consisting of several $p^*$--terms and several $x_{i*}$--terms each.
If such a product contains $\ge 2\ x_{i*}$--terms, it
vanishes because the structure sections pairwise do not intersect.
The product containing no $x_{i*}$--terms vanishes because
$\roman{dim}\ \overline{M}_{0,m+2}=m-1.$ Finally,
there are $m$ products containing one $x_{i*}$--term each.
Using the projection formula
$$
\langle p^*(X)x_{i*}([ \overline{M}_{0,m+2}])\rangle=\langle X\rangle
$$
one sees that each such term equals $-g(m-1)$ (cf. (A.6)).
So $g(m)=-mg(m-1)=(-1)^mm!$ because $g(1)=-1.$

\newpage

\centerline{\bf References}

\bigskip

[AC] E. Arbarello, M. Cornalba. {\it Combinatorial and algebro-geometric
cohomology classes on the moduli spaces of curves.} Preprint, 1994.

\smallskip

[DFI] P. Di Francesco, C. Itzykson. {\it Quantum intersection rings.}
Preprint, 1994.

\smallskip

[D] B. Dubrovin. {\it Geometry of 2D topological field theories.}
Preprint, 1994.

\smallskip

[G] E. Getzler. {\it Operads and moduli spaces of genus zero
Riemann surfaces.} Preprint, 1994.

\smallskip

[GK] E. Getzler, M. Kapranov. {\it Cyclic operads and cyclic
homology.} In ``Geometry, Topology, and Physics for Raoul'',
ed. by B. Mazur, Cambridge MA, 1995.

\smallskip

[Ke] S. Keel. {\it Intersection theory of moduli spaces of stable
$n$--pointed curves of genus zero.} Trans. AMS, 330 (1992), 545--574.

\smallskip

[KM] M. Kontsevich, Yu. Manin. {\it Gromov--Witten classes, quantum
cohomology, and enumerative geometry.} Comm. Math. Phys.,
164:3 (1994), 525--562.

\smallskip

[M] M. Matone. {\it Nonperturbative model of Liouville gravity.}
Preprint hep--th/9402081.

\smallskip

[W] S. Wolpert. {\it On the homology of the moduli spaces of stable
curves.} Ann. of Math., 118 (1983), 491-523.

\smallskip

[Z] P. Zograf. {\it The Weil--Petersson volume of the moduli spaces
of punctured spheres.} Preprint, 1991.

\bigskip

{\it e-mail addresses:}

maxim\@math.berkley.edu

manin\@mpim-bonn.mpg.de

ralph\@mpim-bonn.mpg.de

\enddocument